\documentclass{article}
\usepackage[affil-it]{authblk}
\usepackage[utf8]{inputenc}
\usepackage{lipsum}

\newenvironment{dedication}
  {\clearpage           % we want a new page
   \thispagestyle{empty}% no header and footer
   \vspace*{\stretch{1}}% some space at the top 
   \itshape             % the text is in italics
  }
  {\par % end the paragraph
   \vspace{\stretch{3}} % space at bottom is three times that at the top
   \clearpage           % finish off the page
  }

\usepackage{graphicx}
\usepackage{amssymb, amsmath}
\usepackage{xcolor}
\usepackage[colorinlistoftodos,prependcaption]{todonotes}

\begin{document}
\title{Notes on Low-rank Matrix Factorization}
\author{Yuan Lu, Jie Yang*
 \\ \texttt{\{joyce.yuan.lu,yangjiera\}@gmail.com}.}
\affil{*\ Faculty of EEMCS,\\ 
      Delft University of Technology,\\ 
      Mekelweg 4, 2628 CD Delft, the Netherlands.}
\date{}

\maketitle
\thispagestyle{empty}

\begin{dedication}
Dedicated to Xiao Baobao and Tu Daye. 
\end{dedication}

\section{Introduction}
Low-rank matrix factorization (\textbf{MF}) is an important technique in data science. The key idea of \textbf{MF} is that there exists latent structures in the data, by uncovering which we could obtain a compressed representation of the data. By factorizing an original matrix to low-rank matrices, \textbf{MF} provides a unified method for dimesion reduction, clustering, and matrix completion.

\textbf{MF} has several nice properties: 1) it uncovers latent structures in the data, while addressing the data sparseness problem \cite{koren2008factorization}; 2) it has an elegant probabilistic interpretation \cite{mnih2007probabilistic}; 3) it can be easily extended with domain specific prior knowledge (e.g., homophily in linked data \cite{tang2013exploiting}), thus suitable for various real-world problems; 4) many optimization methods such as (stochastic) gradient-based methods can be applied to find a good solution. 

In this article we review several important variants of \textbf{MF}, including:
\begin{itemize}
    \item Basic \textbf{MF},
    \item Non-negative \textbf{MF},
    \item Orthogonal non-negative \textbf{MF}.
\end{itemize}
As can be seen from their names, non-negative \textbf{MF} and orthogonal non-negative \textbf{MF} are variants of basic \textbf{MF} with non-negativity and/or orthogonality constraints. Such constraints are useful in specific senarios. In the first part of this article, we introduce, for each of these models, the application scenarios, the distinctive properties, and the optimizing method. Note that for the optimizing method, we mainly use the alternative algorithm, as similar to \cite{ding2008nonnegative, tang2013exploiting}. We will derive the updating rules, and prove the correctness and convergence. For reference, matrix operation and optimization can be referred to \cite{brookes2005matrix} and \cite{boyd2009convex} respectively.

%In the context of social computing, \textbf{MF} is widely used in various prediction problems including recommendation, link prediction, sentiment analysis, etc.

By properly adapting \textbf{MF}, we can go beyond the problem of clustering and matrix completion. In the second part of this article, we will extend \textbf{MF} to sparse matrix compeletion, enhance matrix compeletion using various regularization methods, and make use of \textbf{MF} for (semi-)supervised learning by introducing latent space reinforcement and transformation. We will see that \textbf{MF} is not only a useful model but also as a flexible framework that is applicable for various prediction problems. 

\section{Theory}
This section introduces the theory in low-rank matrix factorization. As introduced before, we will go through the following three \textbf{MF} variations: basic \textbf{MF}, non-negative \textbf{MF}, orthogonal non-negative \textbf{MF}.

\subsection{Basic \textbf{MF}}
We start with the basic \textbf{MF} model, formulated as 
\begin{equation}
    \underset{\mathbf{U,V}}{\hbox{min}} \  \|\mathbf{X-UV}^T\|+\mathcal{L}(\mathbf{U,V)},
    \label{BasicMF}
\end{equation}
where $\mathbf{X}\in \mathbb{R}^{m\times n}$ is the data matrix to be approximated, and $\mathbf{U}\in \mathbb{R}^{m\times k}, \mathbf{V}\in \mathbb{R}^{n\times k}$ are two low-dimensional matrices ($k
\ll\hbox{min}(m,n)$). $\mathcal{L}(\mathbf{U,V})$ is a regularization part to avoid overfitting. Regularization is usually necessary in prediction for bias-variance trade-off \cite{hastie2009elements}. 

%This model is widely used in collaborative filtering (CF) towards a recommender system.

\subsubsection{Gradient Descent Optimization}
We instantiate Eq. \ref{BasicMF} as follows
\begin{equation}
    \underset{\mathbf{U,V}}{\hbox{min}} \ \mathcal{O} =  \|\mathbf{X-UV}^T\|_F^2+\alpha\|\mathbf{U}\|_F^2+\beta\|\mathbf{V}\|_F^2.
    \label{BasicMFinit}
\end{equation}

The reason of using Frobenius Norm is that it has a Guassian noise interpretation, and that the objective function can be easily transformed to a matrix trace version:

\begin{equation}
    \underset{\mathbf{U,V}}{\hbox{min}} \ \mathcal{O} = Tr(\mathbf{X}^T\mathbf{X}+\mathbf{VU}^T\mathbf{UV}^T-2\mathbf{X}^T\mathbf{UV}^T)+\alpha Tr(\mathbf{U}^T\mathbf{U})+\beta Tr(\mathbf{V}^T\mathbf{V}).
    \label{BasicMFinitTr}
\end{equation}

Here the matrix calculation rule $\|\mathbf{A}\|_F=\sqrt{Tr(\mathbf{A}^T\mathbf{A})}$  is used in the transformation. Note that trace has many good properties such as  $Tr(\mathbf{A})=Tr(\mathbf{A}^T)$ and $Tr(\mathbf{AB})=Tr(\mathbf{BA})$, which will be used in the following derivations. 

According to trace derivatives $\frac{\partial Tr(\mathbf{AB})}{\partial \mathbf{A}} = \mathbf{B}^T$ and the following rules:
\begin{equation}
\begin{split}
 \frac{\partial Tr(\mathbf{A}^T\mathbf{AB})}{\partial \mathbf{A}} = \mathbf{A}(\mathbf{B}^T+\mathbf{B}), \\
 \frac{\partial Tr(\mathbf{A}\mathbf{A}^T\mathbf{B})}{\partial \mathbf{A}} = (\mathbf{B}^T+\mathbf{B})\mathbf{A}
 \end{split}
\end{equation} (see more in \cite{brookes2005matrix}), we have the following derivatives for $\mathbf{U}$ and $\mathbf{V}$,

\begin{equation}
\begin{split}
\frac{\partial \mathcal{O}}{\partial \mathbf{U}} & = \frac{\partial Tr(\mathbf{VU}^T\mathbf{UV}^T-2\mathbf{X}^T\mathbf{UV}^T)+\alpha Tr(\mathbf{U}^T\mathbf{U})}{\partial \mathbf{U}} \\
 & = \frac{\partial Tr(\mathbf{U}^T\mathbf{UV}^T\mathbf{V}-2\mathbf{UV}^T\mathbf{X}^T)+\alpha Tr(\mathbf{U}^T\mathbf{U})}{\partial \mathbf{U}} \\
 & = 2(\mathbf{UV}^T\mathbf{V} - \mathbf{XV}+\alpha \mathbf{U}), \\
\frac{\partial \mathcal{O}}{\partial \mathbf{V}} & = \frac{\partial Tr(\mathbf{VU}^T\mathbf{UV}^T-2\mathbf{X}^T\mathbf{UV}^T)+\beta Tr(\mathbf{V}^T\mathbf{V})}{\partial \mathbf{V}} \\
 & = \frac{\partial Tr(\mathbf{V}^T\mathbf{VU}^T\mathbf{U}-2\mathbf{V}^T\mathbf{X}^T\mathbf{U})+\beta Tr(\mathbf{V}^T\mathbf{V})}{\partial \mathbf{V}} \\
 & = 2(\mathbf{VU}^T\mathbf{U} - \mathbf{X}^T\mathbf{U}+\beta \mathbf{V}).
\end{split}
\end{equation}

%Both derivations use the commutative law of trace operator. 

Using these two derivatives, we can alternatively update $\mathbf{U}$ and $\mathbf{V}$ in each iteration of gradient descent algorithm. 

Note that the derivation can also be performed elementarily for each entry in matrix $\mathbf{U,V}$ -- this is, in fact, the original definition of matrix calculus. Such element-wise derivation is especially useful in stochastic optimization. We will touch this in a brief discussion of different algorithm schemes next.

\subsubsection{Algorithm Schemes in CF and Others}
For collaborative filtering, usually we take one subset of rated entries in $\mathbf{X}$ as training set, and the rest rated entries as validation set. Detailed algorithm can be found in \cite{takacs2008matrix}. An important implementation strategy is that, for each rated entry in the training set, we update an entire row of $\mathbf{U}$ and an entire column of $\mathbf{V}^T$, as the whole row or column is involved in approximating the rated entry. Same updating mechanism could be applied in stochastic algorithm. 

In the meanwhile, similarly to stochastic algorithm, this type of updating does not fully utilize the data matrix in each updating iteration. The reason is that, not only an entire row of $\mathbf{U}$ ( and a column of $\mathbf{V}^T$) is involved in a single entry in data matrix $\mathbf{X}$, but also that a row of $\mathbf{U}$ (and a column of $\mathbf{V}^T$) influences an entire row (column) of $\mathbf{X}$. Therefore for faster convergence, we recommend to update the matrix $\mathbf{U}$ and $\mathbf{V}$ by fully using data matrix $\mathbf{X}$. 

As the objective function is non-convex caused by the coupling between $\mathbf{U}$ and $\mathbf{V}$, we can choose to alternatively update $\mathbf{U}$ and $\mathbf{V}$ in each iteration as in \cite{ding2008nonnegative, tang2013exploiting}. Detailed algorithm is similar to the one in \cite{tang2013exploiting}. Within any of these matrices, updating should be performed simultaneously as in all gradient-based methods.  Note that, we still need to choose a small learning rate to ensure that the objective function is monotonically decreasing. Interestingly, the alternative optimization scheme is even more suitable for non-negative \textbf{MF} \cite{lee1999learning, lee2001algorithms, ding2006orthogonal, ding2008nonnegative}, as we will see in the following subsections. 

% \vspace{10pt}
% \todo[inline]{Needs comparative validation on large dataset}
% \vspace{10pt}

\subsection{Non-negative \textbf{MF}}
Non-negative \textbf{MF} \cite{lee1999learning} seeks to approximate data matrix $\mathbf{X}$ with low-dimensional matrices $\mathbf{U,V}$ whose entries are all non-negative, i.e., $\mathbf{U,V}\ge \mathbf{0}$. The new problem becomes:

\begin{equation}
\begin{split}
    \underset{\mathbf{U,V}}{\hbox{min}} \ \mathcal{O} &=  \|\mathbf{X-UV}^T\|_F^2+\alpha\|\mathbf{U}\|_F^2+\beta\|\mathbf{V}\|_F^2 \\
    s.t. & \ \mathbf{U}\ge \mathbf{0}, \mathbf{V}\ge \mathbf{0}.
    \label{NMFinit}
\end{split}
\end{equation}

Non-negativity constaint is originated from parts-of-whole interpretation \cite{lee1999learning}. As we can think of, many real-world data are non-negative, such as link strength, favorite strength, etc. Non-negative \textbf{MF} may uncover the important parts, which sometimes can not be achieved by non-constrained \textbf{MF} \cite{lee1999learning}.

Apart from the advantage of uncovering parts, non-negative \textbf{MF} has its own computational advantage: there is a relatively fixed method to find a learning rate larger than common gradient-based methods. To illustrate this, we will first derive the updating rule for Eq. \ref{NMFinit} as an example, then show the general approach for proving the convergence of updating rules derived from the relatively fixed method.

\subsubsection{Updating Rule Derivation}
The basic idea is using KKT complementary slackness conditions to enforce the non-negativity constraint. Based on this, we can directly obtain updating rules. 

The Lagrangian function of Eq. \ref{NMFinit} is
\begin{equation}
    L =  \|\mathbf{X-UV}^T\|_F^2+\alpha\|\mathbf{U}\|_F^2+\beta\|\mathbf{V}\|_F^2 - Tr(\Lambda_1\mathbf{U}^T) - Tr(\Lambda_2\mathbf{V}^T).
    \label{NMFLagrangian}
\end{equation}

We have the following KKT condition,
\begin{equation}
\begin{split}
\Lambda_1 \circ \mathbf{U} = \mathbf{0},  \\
\Lambda_2 \circ \mathbf{V} = \mathbf{0},
\end{split}
\label{NMFKKT1}
\end{equation}
where $\circ$ denotes the Hadamard product. We then have 

\begin{equation}
\begin{split}
\frac{\partial L}{\partial \mathbf{U}} & = \frac{\partial Tr(\mathbf{VU}^T\mathbf{UV}^T-2\mathbf{X}^T\mathbf{UV}^T)+\alpha Tr(\mathbf{U}^T\mathbf{U}) - Tr(\Lambda_1\mathbf{U}^T)}{\partial \mathbf{U}} \\
 & = 2(\mathbf{UV}^T\mathbf{V} - \mathbf{XV}+\alpha \mathbf{U}) - \Lambda_1, \\
\frac{\partial L}{\partial \mathbf{V}} & = \frac{\partial Tr(\mathbf{VU}^T\mathbf{UV}^T-2\mathbf{X}^T\mathbf{UV}^T)+\beta Tr(\mathbf{V}^T\mathbf{V})- Tr(\Lambda_2\mathbf{V}^T)}{\partial \mathbf{V}} \\
 & = 2(\mathbf{VU}^T\mathbf{U} - \mathbf{X}^T\mathbf{U}+\beta \mathbf{V}) - \Lambda_2.
\end{split}
\end{equation}

Let $\frac{\partial L}{\partial \mathbf{U}}=0$ and $\frac{\partial L}{\partial \mathbf{V}}=0$ as another KKT condition, we have

\begin{equation}
\begin{split}
\Lambda_1 & = 2(\mathbf{UV}^T\mathbf{V} - \mathbf{XV}+\alpha \mathbf{U}), \\
\Lambda_2 & = 2(\mathbf{VU}^T\mathbf{U} - \mathbf{X}^T\mathbf{U}+\beta \mathbf{V}).
\end{split}
\label{NMFKKT2}
\end{equation}

Now we combine Eq. \ref{NMFKKT1} and Eq. \ref{NMFKKT2}, we have

\begin{equation}
\begin{split}
(\mathbf{UV}^T\mathbf{V} - \mathbf{XV}+\alpha \mathbf{U})\circ \mathbf{U} = \mathbf{0}, \\
(\mathbf{VU}^T\mathbf{U} - \mathbf{X}^T\mathbf{U}+\beta \mathbf{V})\circ \mathbf{V} = \mathbf{0}.
\end{split}
\label{NMFKKT}
\end{equation}
from which, we have the final updating rules,

\begin{equation}
\begin{split}
\mathbf{U}(i,j) \leftarrow \mathbf{U}(i,j)  \sqrt{\frac{\mathbf{(XV)}(i,j)}{(\mathbf{UV}^T\mathbf{V}+\alpha \mathbf{U})(i,j)}}, \\
\mathbf{V}(i,j) \leftarrow \mathbf{V}(i,j) \sqrt{\frac{(\mathbf{X}^T\mathbf{U})(i,j)}{(\mathbf{VU}^T\mathbf{U}+\beta \mathbf{V})(i,j)}}.
\end{split}
\label{NMFrule}
\end{equation}

Detailed algorithm using these rules is similar to the one in \cite{tang2013exploiting}. We can see that, instead of manually setting small learning rates $\Lambda$'s, Eq. \ref{NMFrule} directly offer updating rules that can usually lead to faster convergence. 

The correctness of these updating rules is straightforward to find out. Taking $\mathbf{U}$ as an example, from Eq. \ref{NMFrule} we have either $\mathbf{U=0}$ or $\mathbf{UV}^T\mathbf{V} - \mathbf{XV}+\alpha \mathbf{U}= \mathbf{0}$, which combined together, exactly equal to Eq. \ref{NMFKKT}.
The convergence, however, is somehow more difficult to be proved. We leave this to the next subsubsection. 
\subsubsection{Proof of Convergence}
We prove the convergence of the updating rules in Eq. \ref{NMFrule} with the standard auxiliary function approach, which is proposed in \cite{lee2001algorithms} and extended in \cite{ding2006orthogonal, ding2008nonnegative}. Our proof is mainly based on \cite{ding2006orthogonal, ding2008nonnegative}, although the objective function Eq. \ref{NMFinit} is slightly different. 

An \textbf{auxiliary function} $G(\mathbf{U}, \mathbf{U}^t)$ of function $L(\mathbf{U})$ is a function that satisfies 
\begin{equation}
G(\mathbf{U}, \mathbf{U}) = L(\mathbf{U}), \  G(\mathbf{U}, \mathbf{U}^t) \ge L(\mathbf{U}).
\label{AuCondition}
\end{equation}
Then, if we take $\mathbf{U}^{t+1}$ such that 
\begin{equation}
\mathbf{U}^{t+1} = \underset{U}{\hbox{arg min}} \  G(\mathbf{U}, \mathbf{U}^t)
\label{AuUpdateVariable},
\end{equation}
we have 
\begin{equation}
L(\mathbf{U}^{t+1})\le G(\mathbf{U}^{t+1}, \mathbf{U}^{t}) \le G(\mathbf{U}^{t}, \mathbf{U}^{t} \le L(\mathbf{U}^{t})).
\end{equation}

This proves that $L(U)$ is monotonically decreasing. 

Turn back to our problem, we need to take two steps using auxiliary function to prove the convergence of updating rules: 1) find an appropriate auxiliary function, and 2) find the global minima of the auxiliary function. As a remark, the auxiliary function approach in principle is similar to Expectation-Maximization approach that is widely used in statistical inference. Now let us complete the proof by taking the above two steps.

\textbf{Step 1 - Finding an appropriate auxiliary function} needs to take advantage of two inequalities,

\begin{align}
z \ge 1 + logz,& \ \forall z>0, \label{PrInequlity1}\\
\sum_{i=1}^m \sum_{j=1}^k \frac{(\mathbf{AS'B})(i,j)\mathbf{S}(i,j)^2}{\mathbf{S'}(i,j)} &\ge Tr(\mathbf{S}^T\mathbf{ASB}), \notag \\
\forall \mathbf{A}\in \mathbb{R}_+^{m\times m}, \mathbf{B}\in \mathbb{R}_+^{k\times k}, &\mathbf{S'}\in \mathbb{R}_+^{m\times k}, \mathbf{S}\in \mathbb{R}_+^{m\times k}.  \label{PrInequlity2}
\end{align}

The proof for Eq. \ref{PrInequlity2} can be found in \cite{ding2006orthogonal} (Proposition 6). 

After removing irrelevant terms, the objective function Eq. \ref{NMFinit} in terms of $\mathbf{U}$ can be written as
\begin{equation}
\begin{split}
& Tr(\mathbf{VU}^T\mathbf{UV}^T-2\mathbf{X}^T\mathbf{UV}^T)+\alpha Tr(\mathbf{U}^T\mathbf{U})\\
= & Tr(\mathbf{U}^T\mathbf{U}\mathbf{V}^T\mathbf{V}-2\mathbf{U}^T\mathbf{X}\mathbf{V})+\alpha Tr(\mathbf{U}^T\mathbf{U}) \\
%= & Tr[-2\mathbf{U}^T\mathbf{X}\mathbf{V}+(\mathbf{X}^T\mathbf{X}+\alpha\mathbf{I})\mathbf{U}^T\mathbf{U}].
\label{PrOriginal}
\end{split}
\end{equation}

We now propose an auxiliary function
\begin{equation}
\begin{split}
G(\mathbf{U}, \mathbf{U}^t) & = -2\sum_{i,j}(\mathbf{XV})(i,j)\mathbf{U}^t(i,j)(1+log\frac{\mathbf{U}(i,j)}{\mathbf{U}^t(i,j)}) \\
 & + \sum_{i,j}\frac{(\mathbf{U}^t\mathbf{V}^T\mathbf{V})(i,j)\mathbf{U}(i,j)^2}{\mathbf{U}^t(i,j)} + \alpha\sum_{i,j}\frac{\mathbf{U}^t(i,j)\mathbf{U}(i,j)^2}{\mathbf{U}^t(i,j)}.
\label{PrAux}
\end{split}
\end{equation}

Combining the two inequalities Eq. \ref{PrInequlity1}, \ref{PrInequlity2}, it is straightforward to see that Eq. \ref{PrAux} is a legal auxiliary function for Eq. \ref{PrOriginal}, i.e., the two conditions in Eq. \ref{AuCondition} are satisfied. Now we procceed to find $\mathbf{U}^{t+1}$ that satisfies condition Eq. \ref{AuUpdateVariable}.

\textbf{Step 2 - Finding $\mathbf{U}^{t+1}$} can be achieved by obtaining the \emph{global minima} of Eq. \ref{PrAux}. First, we have

\begin{equation}
\begin{split}
\frac{\partial G(\mathbf{U}, \mathbf{U}^t)}{\partial \mathbf{U}(i,j)} & = -2(\mathbf{XV})(i,j)\frac{\mathbf{U}^t(i,j)}{\mathbf{U}(i,j)} + 2\frac{(\mathbf{U}^t\mathbf{V}^T\mathbf{V})(i,j)\mathbf{U}(i,j)}{\mathbf{U}^t(i,j)} + 2\alpha\mathbf{U}(i,j).
\end{split}
\end{equation}
Let $\frac{\partial G(\mathbf{U}, \mathbf{U}^t)}{\partial \mathbf{U}(i,j)} = 0$ we have 
\begin{equation}
\begin{split}
(\mathbf{XV})(i,j)\frac{\mathbf{U}^t(i,j)}{\mathbf{U}^{t+1}(i,j)} = (\frac{(\mathbf{U}^t\mathbf{V}^T\mathbf{V})(i,j)}{\mathbf{U}^t(i,j)} + \alpha) \mathbf{U}^{t+1}(i,j),
\end{split}
\end{equation}
from which we directly have
\begin{equation}
\begin{split}
\mathbf{U}^{t+1}(i,j) = \mathbf{U}^t(i,j)\sqrt{\frac{(\mathbf{XV})(i,j)}{(\mathbf{U}^t\mathbf{V}^T\mathbf{V}+\alpha\mathbf{U}^t)(i,j)}},
\end{split}
\label{NMFruleDer}
\end{equation}
which is exactly the updating rule for $\mathbf{U}$ in Eq. \ref{NMFrule}. Similar result can be obtained for $\mathbf{V}$.

\paragraph{General observation} If we go over the entire derivation process, by comparing Eq. \ref{NMFruleDer} and Eq. \ref{NMFKKT}, we can observe that the only thing that matters for the final updating rules is the signs of the terms in Eq. \ref{NMFKKT}.

\subsection{Orthogonal Non-negative \textbf{MF}}
Orthogonality is another important constraint to \textbf{MF}. First of all, we formulate the problem as

\begin{equation}
\begin{split}
    \underset{\mathbf{U,V}}{\hbox{min}} \ \mathcal{O} &=  \|\mathbf{X-UV}^T\|_F^2 \\
    s.t. & \ \mathbf{U}, \mathbf{V} \ge \mathbf{0}, \mathbf{U}^T\mathbf{U} = \mathbf{I}, \mathbf{V}^T\mathbf{V} = I.
    \label{ONMFinit}
\end{split}
\end{equation}
Note that here we do not add regularization due to the orthogonality constraint. 

It is proved in \cite{ding2005equivalence, ding2006orthogonal} (\cite{ding2006orthogonal} gives more mature proof) that this problem is equivalent to K-means clustering: $\mathbf{V'}$ is an indication matrix with $\mathbf{V'}(i,j) = 0$ if $\mathbf{x}_i$ belongs to the $j^{th}$ $(1\leq j \leq k)$ cluster. Here $\mathbf{V}=\mathbf{V'}(\mathbf{V'}^T\mathbf{V'})^{-1/2}$, i.e., $\mathbf{V}$ is a normalized version of $\mathbf{V'}$: $\mathbf{V'}$ is a constant scaling of corresponding row of $\mathbf{V}$, and $\|\mathbf{V}(:,j)\|_2^2 = 1$.

\subsubsection{3-factor \textbf{MF} vs. 2-factor \textbf{MF}}
We call Eq. \ref{ONMFinit}  1-sided 2-factor orthogonal non-negative \textbf{MF}, as only one factorized matrix needs to be orthogonal, and there are in total two factorized matrices.  It is  recommended that, to simultaneously cluster rows and columns in $\mathbf{X}$, we need 3-factor bi-orthogonal non-negative \textbf{MF}, i.e.,  both $\mathbf{U}$ and $\mathbf{V}$ being orthogonal:

\begin{equation}
\begin{split}
    \underset{\mathbf{U,H,V}}{\hbox{min}} \ \mathcal{O} &=  \|\mathbf{X-UHV}^T\|_F^2 \\
    s.t. & \ \mathbf{U}, \mathbf{H}, \mathbf{V} \ge \mathbf{0}, \mathbf{U}^T\mathbf{U} = \mathbf{I}, \mathbf{V}^T\mathbf{V} = I.
    \label{3ONMFinit}
\end{split}
\end{equation}

It is proved that, compared to 3-factor bi-orthogonal non-negative \textbf{MF}, 2-factor bi-orthogonal non-negative \textbf{MF} is too restrictive, and will lead to poor approximation \cite{ding2006orthogonal}. 

3-factor bi-orthogonal non-negative \textbf{MF} is useful in document-word clustering \cite{ding2006orthogonal}, outperforming K-means (i.e., 1-sided 2-factor orthogonal non-negative \textbf{MF}). It has been applied for tasks such as sentiment analysis \cite{hu2013unsupervised}.

\subsubsection{Updating Rule Derivation}
We now derive updating rules for Eq. \ref{3ONMFinit}, as we did before for non-negative \textbf{MF}. 

The Lagrangian function for Eq. \ref{3ONMFinit} is

\begin{equation}
\begin{split}
   L = & \|\mathbf{X-UHV}^T\|_F^2 - Tr(\Lambda_U\mathbf{U}^T) - Tr(\Lambda_H\mathbf{H}^T) - Tr(\Lambda_V\mathbf{V}^T) \\+& Tr(\Gamma_U(\mathbf{U}^T\mathbf{U} - \mathbf{I})) + Tr(\Gamma_V(\mathbf{V}^T\mathbf{V} - \mathbf{I}))
    \label{3ONMFinitLag}
\end{split}
\end{equation}

We then compute the updating rules for $\mathbf{H,U,V}$ sequentially.

\textbf{Computation of} $\mathbf{H}$
\begin{equation}
\begin{split}
\frac{\partial L}{\partial \mathbf{H}} & = \frac{\partial Tr(\mathbf{VH}^T\mathbf{U}^T\mathbf{UHV}^T-2\mathbf{XVH}^T\mathbf{U}^T)-Tr(\Lambda_H\mathbf{H}^T)}{\partial \mathbf{H}}\\
& = 2\mathbf{U}^T\mathbf{U}\mathbf{H}\mathbf{V}^T\mathbf{V}-2\mathbf{U}^T\mathbf{XV}-\Lambda_H,
\end{split}
\end{equation}

We have the following KKT conditions,
\begin{equation}
\begin{split}
\frac{\partial L}{\partial \mathbf{H}} & = \mathbf{0}\\
\Lambda_H\circ \mathbf{H} & = \mathbf{0}.
\end{split}
\end{equation}

Combining the above three equations, we have 
\begin{equation}
\begin{split}
(\mathbf{U}^T\mathbf{U}\mathbf{H}\mathbf{V}^T\mathbf{V}-\mathbf{U}^T\mathbf{XV})\circ \mathbf{H} = \mathbf{0}.
\end{split}
\end{equation}

Therefore we have the following updating rule for $\mathbf{H}$,
\begin{equation}
\begin{split}
\mathbf{H}(i,j) \leftarrow \mathbf{H}(i,j)  \sqrt{\frac{(\mathbf{U}^T\mathbf{XV})(i,j)}{(\mathbf{U}^T\mathbf{U}\mathbf{H}\mathbf{V}^T\mathbf{V})(i,j)}}.
\end{split}
\end{equation}

Note that $\mathbf{U}^T\mathbf{U}\ne \mathbf{I}$ during the optimizing process.

\textbf{Computation of} $\mathbf{U,V}$ 

Due to the orthogonality constraint, obtaining the updating rules for $\mathbf{U, V}$ needs to eliminate both $\Lambda$ and $\Gamma$ in the final updating rules. This will need the following equality,

\begin{equation}
\mathbf{U}^T\Lambda_U = \mathbf{0} \Leftarrow \Lambda_U\circ\mathbf{U} = \mathbf{0}
\label{ONMFequality}
\end{equation}
The latter will automatically be satisifed according to KKT conditions as we will see below. 

\begin{equation}
\begin{split}
\frac{\partial L}{\partial \mathbf{U}} & = \frac{\partial Tr(\mathbf{VH}^T\mathbf{U}^T\mathbf{UHV}^T-2\mathbf{XVH}^T\mathbf{U}^T)-Tr(\Lambda_U\mathbf{U}^T) + Tr(\Gamma_U(\mathbf{U}^T\mathbf{U} - \mathbf{I}))}{\partial \mathbf{U}}\\
& = 2\mathbf{UHV}^T\mathbf{V}\mathbf{H}^T-2\mathbf{X}\mathbf{VH}^T-\Lambda_U+2\mathbf{U}\Gamma_U,
\end{split}
\end{equation}

We have the following KKT conditions,
\begin{equation}
\begin{split}
\frac{\partial L}{\partial \mathbf{U}} & = \mathbf{0}\\
\Lambda_U\circ \mathbf{U} & = \mathbf{0}.
\end{split}
\end{equation}

Combining the above three equations we have
\begin{equation}
\begin{split}
(\mathbf{UHV}^T\mathbf{V}\mathbf{H}^T-\mathbf{X}\mathbf{VH}^T +\mathbf{U}\Gamma_U)\circ \mathbf{U} = \mathbf{0} 
\label{ONMFlastbuttwo}
\end{split}
\end{equation}
and
\begin{equation}
\begin{split}
\Gamma_U = \mathbf{U}^T\mathbf{X}\mathbf{VH}^T - \mathbf{HV}^T\mathbf{V}\mathbf{H}^T.
\label{ONMFlastbutone}
\end{split}
\end{equation}

Note that here we can have $\mathbf{U}^T\mathbf{U}=\mathbf{I}$ as we only want an expression for $\Gamma_U$. Further note that for $\Lambda$ we have the constraint $\Lambda>\mathbf{0}$ (according to KKT condition) while for $\Gamma$ we do not have such constraint. Therefore we need to split $\Gamma$ into two parts,
\begin{equation}
\begin{split}
\Gamma_U &= \Gamma_U^+ - \Gamma_U^- \\
\Gamma_U^+ &= (|\Gamma_U|+\Gamma_U)/2 \\
\Gamma_U^- &= (|\Gamma_U|-\Gamma_U)/2.
\end{split}
\label{splitruleU}
\end{equation}

Using this division we rewrite Eq. \ref{ONMFlastbuttwo}, we then have

\begin{equation}
\begin{split}
(\mathbf{UHV}^T\mathbf{V}\mathbf{H}^T-\mathbf{X}\mathbf{VH}^T +\mathbf{U}\Gamma_U^+ - \mathbf{U}\Gamma_U^-)\circ \Lambda_U = \mathbf{0}.
\end{split}
\end{equation}

Therefore the final updating rule for $\mathbf{U}$ is
\begin{equation}
\begin{split}
\mathbf{U}(i,j) \leftarrow \mathbf{U}(i,j)  \sqrt{\frac{(\mathbf{X}\mathbf{VH}^T+\mathbf{U}\Gamma_U^-)(i,j)}{(\mathbf{UHV}^T\mathbf{V}\mathbf{H}^T+\mathbf{U}\Gamma_U^+)(i,j)}}.
\end{split}
\end{equation}
where $\Gamma_U^+$ and $\Gamma_U^-$ is defined in Eq. \ref{ONMFlastbutone} and \ref{splitruleU}.

If we go over the same process again for $\mathbf{V}$, we have the following updating rules,
Therefore the final updating rule for $\mathbf{U}$ is
\begin{equation}
\begin{split}
\mathbf{V}(i,j) \leftarrow \mathbf{V}(i,j)  \sqrt{\frac{(\mathbf{X}^T\mathbf{UH}+\mathbf{V}\Gamma_V^-)(i,j)}{(\mathbf{VH}^T\mathbf{U}^T\mathbf{UH}+\mathbf{V}\Gamma_V^+)(i,j)}}.
\end{split}
\end{equation}
where $\Gamma_V^+,\Gamma_V^-$ are defined similarly as in Eq. \ref{splitruleU} (replace $\mathbf{U}$ with $\mathbf{V}$), and $\Gamma_V$ is defined as  
\begin{equation}
\Gamma_V =\mathbf{V}^T\mathbf{X}^T\mathbf{UH}+\mathbf{H}^T\mathbf{U}^T\mathbf{UH}.
\end{equation}

\paragraph{Choice of 2/3-factor \textbf{MF}} How do we choose between 2-factor or 3-factor \textbf{MF} in real-world applications? A general principle is that: if we only need to place regularizations on one latent matrix, i.e. either $\mathbf{U}$ or $\mathbf{V}$, then we can use 2-factor \textbf{MF}; if both $\mathbf{U}$ and $\mathbf{V}$ are to be regularized, either explictly or implictly, 3-factor \textbf{MF} might be a better choice.

\section{Adapatations and Applications}
\textbf{MF} has been used for a wide range of applications in social computing, including collaborative filtering (CF), link prediction (LP), sentiment analysis, etc. It can not only provide as a single model for matrix completeion or clutering, but also as a framework for solving almost all categories of prediction problems.

In this part we will extend \textbf{MF} to highly sparse cases. For the cases in which we have additional data, e.g. link data between users (in CF, or addtional links in LP) or description data of users and items, we can incorporate different regularization techniques to enhace the matrix completion performance. Moreover, by properly manipulating latent factors derived from \textbf{MF}, we can adapt \textbf{MF} to (un-/semi-)supervised learning.

\subsection{Sparse Matrix Completion}
Here we address the problem of using \textbf{MF} for collborative filtering, link prediction and clustering. We start with a basic assumption, which makes the previously introduced models unsuitable. This basic assumption is: high portion of the data is missing, i.e. data matrix is incomplete. Such assumption is very common in real-world cases \cite{koren2009matrix}. 

The problem is solved by modeling directly the observed data. Eq. \ref{BasicMF} is modified as follows:
\begin{equation}
    \underset{\mathbf{U,V}}{\hbox{min}} \ \mathcal{O} =  \|\mathbf{O\circ (X-UV}^T)\|_F^2+\alpha\|\mathbf{U}\|_F^2+\beta\|\mathbf{V}\|_F^2,
    \label{BasicMFinitCF}
\end{equation} in which $O$ poses constraints on only these observed data entries, i.e. $\mathbf{O}(i,j)=1$ if entry $(i,j)$ is observed, and $\mathbf{O}(i,j)=0$ otherwise. 

In this case, the objective function is transformed as follows:
\begin{equation}
\begin{split}
    \underset{\mathbf{U,V}}{\hbox{min}} \ \mathcal{O} &= Tr((\mathbf{O}^T\circ \mathbf{X}^T)(\mathbf{O\circ X})+(\mathbf{O}^T\circ \mathbf{VU}^T)(\mathbf{O\circ UV}^T)\\&-2(\mathbf{O}^T\circ \mathbf{X}^T)(\mathbf{O\circ UV}^T))+ \alpha Tr(\mathbf{U}^T\mathbf{U})+\beta Tr(\mathbf{V}^T\mathbf{V}).
    \label{BasicMFinitTrCF}
\end{split}
\end{equation}

And the gradients become:
\begin{equation}
\begin{split}
\frac{\partial \mathcal{O}}{\partial \mathbf{U}} & = \frac{\partial Tr((\mathbf{O}^T\circ \mathbf{VU}^T)(\mathbf{O\circ UV}^T)-2(\mathbf{O}^T\circ \mathbf{X}^T)(\mathbf{O\circ UV}^T))+\alpha Tr(\mathbf{U}^T\mathbf{U})}{\partial \mathbf{U}} \\
 & = \frac{\partial Tr(\mathbf{U}^T(\mathbf{O\circ O\circ UV}^T)\mathbf{V}-2(\mathbf{O}^T\circ \mathbf{O}^T\circ \mathbf{X}^T)\mathbf{UV}^T)+\alpha Tr(\mathbf{U}^T\mathbf{U})}{\partial \mathbf{U}} \\
 & = 2((\mathbf{O\circ O\circ UV}^T)\mathbf{V} - (\mathbf{O}\circ \mathbf{O}\circ \mathbf{X})\mathbf{V}+\alpha \mathbf{U}), \\
\frac{\partial \mathcal{O}}{\partial \mathbf{V}} &  = 2((\mathbf{O}^T\circ \mathbf{O}^T\circ  \mathbf{VU}^T)\mathbf{U} - (\mathbf{O}^T\circ \mathbf{O}^T \circ \mathbf{X}^T)\mathbf{U}+\beta \mathbf{V}).
\end{split}
\end{equation}

In the derivation above we use the following rule of Hadamard product:
\begin{equation}
Tr((\mathbf{O}^T \circ \mathbf{A}^T)(\mathbf{O}\circ \mathbf{A})) = Tr(\mathbf{A}^T(\mathbf{O}\circ \mathbf{O}\circ \mathbf{A})).
\end{equation}

The upodating rules for non-negative \textbf{MF} and orthogonal non-negative \textbf{MF} is straightforward: the methods of getting $\Lambda, \Gamma$ are exactly the same as what we did in Theory Section. For updating rules of non-negative \textbf{MF} and orthogonal non-negative \text{MF}, the reader can refer to \cite{gao2013exploring} and \cite{gao2015content}, respectively. 

\subsubsection{Calculating Memory Occupation}
Note that the updating rules above are again purely matrix-wise -- this is to be consistent with the style of this article. In matrix completion, however, sometimes the size of the data matrix is bigger than memory size, making stochasitc gradient descent algorithm more suitable than the matrix-wise method.  

The question here is, how do we calculate the size of a matrix to see if it fits to memory. Here is a easy way to make such a calculation. Assume we have a $10K\times 10K$ matrix, with each entry allocated a 32bit float (e.g. float32 in python), then the memory allocation for the whole matrix can be roughtly calculated as 
\[(10^4\times 10^4\times 4)/10^6 = 400M.\]

So for a computer with 4G memory, we can fit a matrix $100K\times 10K$ matrix into memory. For a computer with 32G memory, we can fit a matrix of size $100K\times 80K$ ($10\times8\times 400M=32G$).

\subsection{Enhanced Matrix Completion}
We looked at \textbf{MF} with different \textit{constraints}, e.g. non-negativity and orthogality, and one type of \textit{regularization} which prevents the entries in low-rank matrices being too large. This subsection considers other kinds of \textit{regularization} when external data source becomes avaiable, i.e. goes beyond the data matrix $\mathbf{X}$.  Usually this is the real-world case, since most social media data contains rich data sources. 

In this subsection we consider two types of regularization with corresponding addtional data:
\begin{enumerate}
\item \textbf{self-regularization} when we have additional linked data between users (in CF, or addtional link type in LP);
\item \textbf{2-sided regularization} when we have description data of users and items.
\end{enumerate}

We further point to two publications \cite{tang2013exploiting} and  \cite{gao2015content}, to demonstrate the above two types of regularization, respectively.

\subsubsection{Enhancing Matrix Completion with Self-regularization}
By self-reguarization, we refer to the regularization of rows in low-rank matrix $\mathbf{U}$ or $\mathbf{V}$. Assume now we are dealing with a LP problem, in which we would like to predict if a user trust another -- trust relation are common in review sites like Epinions. Usually there exist another type of links between users, i.e. social relation.  Can we use social relation to boost the performance of trust relation prediction? This is exactly the research question proposed in \cite{tang2013exploiting}. 

It turns out the answer is yes -- as expected, users with social relation tend to share similar preferences. The basic idea to incorporate this into trust prediction is by adding the  regularization term Eq.~\ref{Laplacian} into the general \textbf{MF} framework. In Eq.~\ref{Laplacian}, $\xi$ is the entries in the additional link matrix $\mathcal{Z}$ and $\mathbf{D}$ is the diagonal matrix with $\mathbf{D}(i,i)=\mathcal{Z}_{j=1}^m(j,i)$, thus $\mathcal{L}$ is the Laplacian matrix of $\mathbf{D}$. It is interesting that, using trace operator, the regularization Eq.~\ref{Laplacian} become such simple.

Social relation is common in social computing, the similarity in people with social relation has a specific name in social theory - `homophily', making this type of regularization applicable to a lot of social computing scenarios. If we generalize a bit, we may assume that many linked objects, not necessarily web users, have similarities, in terms of their entries of data matrix $\mathbf{X}$ that we would like to predict. For instance, while predicting the sentiment of articles, we may assume that articles authored by the same users tend to express similar sentiment, e.g. political reviewers expressing negative sentiment in their news reviewing articles. We will see that this type of regularization is used in a sentiment analysis paper \cite{hu2013unsupervised}, which we will analyze later.

\begin{equation}
\begin{split}
	&\frac{1}{2}\sum_{i=1}^m\sum_{j=1}^m \xi(i,j) \|\mathbf{U}(i,:)-\mathbf{U}(j,:)\|_2^2 \\
=	&\frac{1}{2}\sum_{i=1}^m\sum_{j=1}^m\sum_{d=1}^k \xi(i,j) (\mathbf{U}(i,k)-\mathbf{U}(j,k))^2 \\
=	&\frac{1}{2}\sum_{i=1}^m\sum_{j=1}^m\sum_{d=1}^k \xi(i,j) (\mathbf{U}^2(i,k)-2\mathbf{U}(i,k)\mathbf{U}(j,k)+\mathbf{U}^2(j,k)) \\
=	&\sum_{i=1}^m\sum_{j=1}^m\sum_{d=1}^k \xi(i,j)\mathbf{U}^2(i,k) - \sum_{i=1}^m\sum_{j=1}^m\sum_{d=1}^k \xi(i,j)\mathbf{U}(i,k)\mathbf{U}(j,k) \\
=	& \sum_{d=1}^k \mathbf{U}^T(:,k)(\mathbf{D}-\mathcal{Z})\mathbf{U}(:,k) \\
=	&Tr(\mathbf{U}^T\mathcal{L}\mathbf{U})
\label{Laplacian}
\end{split}
\end{equation}  

\paragraph{Regularization and Sparseness} More regularization sometimes can conquer the data sparsity problem, to some extent. On the other hand, modelling the error only on observed data entries, as what $\mathbf{O}$ does in previous subsection, could be also very effective.

\subsubsection{Enhancing Matrix Completion with 2-sided regularization}
Here we consider placing regularization on both $\mathbf{U}$ and $\mathbf{V}$ together, which we call 2-sided regularization. 

Before we start, we review orthogonal non-negative \text{MF} a bit. Orthogonality constraint in orthogonal non-negative \text{MF} is similar to a 2-sided regularization:
\[Tr(\Gamma_U^T(\mathbf{U}^T\mathbf{U} - \mathbf{I})), Tr(\Gamma_V^T(\mathbf{V}^T\mathbf{V} - \mathbf{I}))\] are two equality constraints over low-rank matrices. Such equality needs to be strictily satisfied. \textbf{Regularization, differing from constraints, however can be viewed as a soft type of constraints}: it only needs to be satisfied to some extend, while constraints need to be strictly satisified. This is the reason why we consider non-negativity and orthogonality constraints, while call homophily regularization. 

Now let us turn our attention back to 2-sided regularization, basing the example from \cite{gao2015content}, which considers POI recommendation in location-based social network (LBSN). The first data we have is a check-in data $\mathbf{X}$ that encodes the interaction between users and POI's.  We are further given some desription data $\mathbf{A}$ of user interest, and $\mathbf{B}$ of POI property, both in the form of word vectors. Question here is, how do we make use of $\mathbf{A}$ and $\mathbf{B}$ to enhance the matrix completion problem for interacting matrix $\mathbf{X}$?

Since we are coping with 2-sided regularization, we use 3-factor \text{MF}:

\begin{equation}
\begin{split}
    \underset{\mathbf{U,H,V}}{\hbox{min}} \ \mathcal{O} &=  \|\mathbf{X-UHV}^T\|_F^2 - Tr(\Lambda_U\mathbf{U}^T) - Tr(\Lambda_H\mathbf{H}^T) + R's.
\end{split}
\end{equation} The only thing here is, how to add the 2-sided regularization terms $R$'s, as we did for orthogonality constraints.

To utilize $\mathbf{A}$ and $\mathbf{B}$, we assume that there are some connections between them, such that they can be used to regularize $\mathbf{U}$ and $\mathbf{V}$. In the context of LBSN, we may assume that $\mathbf{A}$ and $\mathbf{B}$ have similar vocabulary, in which the words have similar latent space. Therefore we can approximate $\mathbf{A}$ and $\mathbf{B}$ with 2-factor \textbf{MF}:
\begin{equation}
\mathbf{A} \approx \mathbf{U}\mathbf{G}^T, \mathbf{B}  \approx \mathbf{V}\mathbf{G^*}^T
\end{equation}
with connection 
\begin{equation}
\|\mathbf{G}-\mathbf{G^*}\|_1 \approx 0.
\label{connect}
\end{equation}

Eq.~\ref{connect} is important since it really connect $\mathbf{U}$ with $\mathbf{V}$, forming a 2-sided regularization. The final objective function now becomes:

\begin{equation}
\begin{split}
    \underset{\mathbf{U,H,V}}{\hbox{min}} \ \mathcal{O} &=  \|\mathbf{X-UHV}^T\|_F^2 - Tr(\Lambda_U\mathbf{U}^T) - Tr(\Lambda_H\mathbf{H}^T) \\& + \lambda_A \|\mathbf{A} - \mathbf{U}\mathbf{G}^T\|_F^2 + \lambda_B \|\mathbf{B} - \mathbf{U}\mathbf{G^*}^T\|_F^2  + \delta \|\mathbf{G}-\mathbf{G^*}\|_1 \\& + \alpha (\|\mathbf{U}\|_F^2+\|\mathbf{V}\|_F^2+\|\mathbf{H}\|_F^2+\|\mathbf{G}\|_F^2).
\end{split}
\end{equation} 

The last line is to regularize in approximating $\mathbf{A,B}$; note that since here we use regularization, instead of constraints as in non-negative orthogonal \textbf{MF}, we can add regualrization to $\mathbf{U,V,H}$.

\paragraph{Factorization vs. Regularization} We remark here that the idea of co-factoring two matrices ($\mathbf{X,A}$) with shared factors ($\mathbf{U}$) originates from collective matrix facterization \cite{singh2008relational}, which has many applications in CF \cite{shi2014collaborative}. A interesting comparative study between collective facterization and self-regularization can be found in \cite{yuan2011factorization}.

\subsection{From Clustering to (Un-/Semi-)supervised Learning}
Although different types of extra data sources can be used in enhanced \textbf{MF}, the purpose  so far to remains be matrix completion. This subsection, however, considers other types of machine learning problems, i.e. (un-/semi-)supervised learning. The essential assumption of using \textbf{MF} for (un-/semi-)supervised learning is that the latent row(column) is or can be predictable for some dependent variables. 

To make use of the predictability, we need mechanisms to connect the latent vectors to responses. Following are the two mechanisms:
\begin{enumerate}
\item \textbf{enforcement} directly enforce the latent space to be the response space;
\item \textbf{transformation} transform the latent space to response space. This is similar as what people do in machine learning.
\end{enumerate}

We point to publications \cite{hu2013unsupervised} and  \cite{gao2014modeling} for the demonstration of the above two methods, respectively. 

\subsubsection{Enforcing Latent Factor to be Response}
In previous regularizations, we do not force the latent space to be interpretable space. For instance, in the 2-sided regularization, we do not specify the meaning of $\mathbf{U}$ that is used in both $\mathbf{X}$ and $\mathbf{A}$ factorization. However, (un-/semi-)supervised learning requires the latent space to be interpretable. The method, still, is regularization. 

\cite{hu2013unsupervised} deals with the problem of sentiment analysis, for which the authors use 3-factor non-negative orthogonal \textbf{MF}. The input is a post-word matrix $\mathbf{X}$. In addition, we are given emotion indication in some of the posts. ``The key idea of modeling post-level emotion indication is to make the sentiment polarity of a post as close as possible to the emotion indication of the post.'', formulated as 
\[\mathbf{G}^u \|\mathbf{U-U_0}\|_F^2,\]
in which $\mathbf{U}\in \mathbb{R}^{m\times 2}$ is the post-sentiment matrix, i.e. $\mathbf{U}(i,:)=(1,0)$ representing that the  $i$th post has a positive sentiment, and $\mathbf{U_0}\in \mathbb{R}^{m\times 2}$ is the post-emotion indication matrix, i.e. $\mathbf{U_0}(i,:)=(1,0)$ meaning the  $i$th post contains positive emotion indication. Similar regularization is applied to $\mathbf{V}$ as well. 

Such an idea is quite simple, however it explictly poses a notable question: is it computationally feasible that we strictly enforce the $\mathbf{U,V}$ to any pre-defined space, i.e. sentiment space in this case. Based on Proposition 1 in \cite{ding2006orthogonal}, we know that the answer is no. However, as we see in this sentiment analysis work \cite{hu2013unsupervised}, \textbf{regularization is always possible}! 

In fact, the enforcement regularization that we see in this work is the most constrained regularization: it is 2-sided regularization for both $\mathbf{U,V}$, and it is enforcement without any transformation coefficients. We will see next how to regularize for supervised learning by tranformation. 

\subsubsection{Transforming Latent Factor to Response}
As we pointed out, the essential idea of supervised learning is to transform the latent variables to some response variable. To see this, we study an example that exploit matrix factorization to boost (sparse) regression.
%\setcounter{equation}{0}
%\newpage
%\textbf{Derivation of Updating Rules for MFLR}

Here we solve the following optimization problem:
\begin{equation}\small
\begin{split}
\underset{\mathbf{U,V}\ge \mathbf{0}}{\hbox{min}}\ \  & \|\mathbf{X- UV}^T\|_F^2 + \lambda\|\mathbf{O}\odot (\mathbf{UW}^T - \mathbf{Y})\|_F^2  \\
 & +  \lambda_X(\|\mathbf{U}\|_F^2+\|\mathbf{V}\|_F^2) + \lambda_Y\|\mathbf{W}\|_1.
\end{split}
\end{equation}

Optimizing the objective function accomplishes two goals simultaneously: 1) learning the latent  factors; and, 2) predicting the dependent variables based on the learnt latent factors. As the learning of $\mathbf{U}$ is guided by the prediction of $\mathbf{Y}$ (proved later), the learned latent factors can be more predictive in the regression. Note that the parameter $\lambda$ controls the relative importance between matrix factorization and regression – a larger  $\lambda$ indicates that the regression should dominate.

$\mathbf{O}$ is a mask vector with the first $n_{train}$ – the size of training set – entries equal to 1, and the other  $n_{test}$ – the size of test set – entries equal to 0. Correspondingly, $\mathbf{X}$ contains both the training data in the first $n_{train}$ rows and the test data, in the remaining $n_{test}$ rows. $\mathbf{Y}$ is also composed of two parts, the first $n_{train}$  entries being the complexity values of the training tasks; the other entries can be any values, as they are not involved in model learning, which is controlled by the 0’s in $\mathbf{O}$.
%We use an alternative algorithm, as similar to \cite{ding2008nonnegative}, to solve this problem. The algorithm to iteratively update $\mathbf{U,V,W}$ until the objective function converges. Since the objective function is non-convex, it alternatively updates $\mathbf{U,V,W}$ in each iteration. %We then derive the updating rules for the variables involved in this problem.

The Lagrangian function of the objective function is
\begin{equation}\small
\begin{split}
L&=\|\mathbf{X- UV}^T\|_F^2 + \lambda\|\mathbf{O}\odot (\mathbf{UW}^T - Y)\|_F^2  
\\ &+ \lambda_X(\|\mathbf{U}\|_F^2+\|\mathbf{V}\|_F^2)
+ \lambda_Y\|\mathbf{W}\|_1 - Tr(\Lambda_U\mathbf{U}^T) - Tr(\Lambda_V\mathbf{V}^T).
\end{split}
\end{equation}

The derivative of $\mathbf{U}$ is:
\begin{equation}\small
\begin{split}
\frac{\partial L}{\partial \mathbf{U}}  
&= \frac{\partial Tr(\mathbf{VU}^T\mathbf{UV}^T-2\mathbf{X}^T\mathbf{UV}^T) }{\partial \mathbf{U}} \\
& + \frac{\partial \lambda Tr((\mathbf{O}^T\odot \mathbf{WU}^T)(\mathbf{O}\odot \mathbf{UW}^T)-2(\mathbf{O}^T\odot \mathbf{Y}^T)(\mathbf{O}\odot \mathbf{UW}^T))}{\partial \mathbf{U}} \\
& + \frac{\partial \lambda_X Tr(\mathbf{U}^T\mathbf{U}) - Tr(\Lambda_U\mathbf{U}^T)}{\partial \mathbf{U}} \\
& = 2(\mathbf{UV}^T\mathbf{V} - \mathbf{XV}+ \lambda(\mathbf{O}\odot (\mathbf{UW}^T))\mathbf{W} - \lambda(\mathbf{O}\odot \mathbf{Y})\mathbf{W} + \lambda_X \mathbf{U}) \\
& - \Lambda_U.
  \end{split}
\end{equation}
  
The derivative of $\mathbf{V}$ is:
\begin{equation}\small
\begin{split}
\frac{\partial L}{\partial \mathbf{V}}  &= \frac{\partial Tr(\mathbf{VU}^T\mathbf{UV}^T-2\mathbf{X}^T\mathbf{UV}^T)+\lambda_X Tr(\mathbf{V}^T\mathbf{V})- Tr(\Lambda_V\mathbf{V}^T)}{\partial \mathbf{V}} \\
 &= 2(\mathbf{VU}^T\mathbf{U} - \mathbf{X}^T\mathbf{U}+\lambda_X \mathbf{V}) - \Lambda_V.
  \end{split}
\end{equation}
 
 Note that for $\mathbf{W}$ the problem becomes a classic Lasso problem, we can update it use standard algorithm such as LARS.
 
 According to the KKT conditions:
\begin{equation}\small
\begin{split}\frac{\partial{L}}{\partial{\mathbf{U}}}=0 &,  \frac{\partial{L}}{\partial{\mathbf{V}}}=0, \\ 
\Lambda_U\odot \mathbf{U}=0 &, \Lambda_V\odot \mathbf{V}=0.
  \end{split}
\end{equation}

We have
\begin{equation}\small
\begin{split}
 & (\mathbf{UV}^T\mathbf{V} - \mathbf{XV}+  \lambda(\mathbf{O}\odot (\mathbf{UW}^T))\mathbf{W}  \\
 & - \lambda(\mathbf{O}\odot \mathbf{Y})\mathbf{W} + \lambda_X \mathbf{U})\odot \mathbf{U} = 0,\\
 & (\mathbf{VU}^T\mathbf{U} - \mathbf{X}^T\mathbf{U}+\lambda_X \mathbf{V})\odot \mathbf{V} = 0. \\
  \end{split}
\end{equation}

It leads to the following updating rules for $\mathbf{U,V}$:
\begin{equation}\small
\begin{split}
\mathbf{U}(i,j) &\leftarrow \mathbf{U}(i,j)  \sqrt{\frac{(\mathbf{XV} +  \lambda(\mathbf{O}\odot \mathbf{Y})\mathbf{W})(i,j)}{(\mathbf{UV}^T\mathbf{V}+ \lambda(\mathbf{O}\odot  (\mathbf{UW}^T))\mathbf{W} + \lambda_X \mathbf{U})(i,j)}}, \\
\mathbf{V}(i,j)& \leftarrow \mathbf{V}(i,j) \sqrt{\frac{(\mathbf{X}^T\mathbf{U})(i,j)}{(\mathbf{VU}^T\mathbf{U}+\lambda_X \mathbf{V})(i,j)}} .
 \end{split}
\end{equation}

\subsubsection{A Comprehensive Model}
Here we review an application of \cite{gao2014modeling} that integrates the methods of enforcement and transformation. In this application, we would like to model a user' attitude towards some controversial topic, reflected by his opinion, sentiment and retweeting action. We are given a retweeting matrix $\mathbf{X}$ representing users' retweeting action to some tweets, and we would like to predict users' opinion $\mathbf{O}$ and sentiment $\mathbf{P}$, and the task is to predict these three variables given the user feature $\mathbf{F}$. 

We first introduce how the model is built in \cite{gao2014modeling}, then discuss other alternatives. To train such a model, the authors propose the following model

\begin{equation}
\begin{split}
    \underset{\mathbf{W,V}}{\hbox{min}} \ \mathcal{O}  & =  \|\mathbf{X}-(\mathbf{FW}^T)\mathbf{V}^T\|_F^2 + \lambda_1 \|\mathbf{FW}^T - \mathbf{O}\|_F^2 + \lambda_2 \|(\mathbf{FW}^T)\mathbf{S} - \mathbf{P}\|_F^2\\& +\lambda_3\|W\|_1+\alpha\|\mathbf{W}\|_F^2+\beta\|\mathbf{V}\|_F^2 +\gamma\|\mathbf{S}\|_F^2 - Tr(\Lambda_1\mathbf{U}^T) - Tr(\Lambda_2\mathbf{V}^T),
    \label{Transform}
\end{split}
\end{equation} in which $\lambda_1 \|\mathbf{FW}^T - \mathbf{O}\|_F^2$ and $\lambda_3\|W\|_1$ models opinion from the user feature by bringing in the classical linear regression. We can see that modeling the sentiment is also straightforward: $\lambda_2 \|\mathbf{FW}^T\mathbf{S} - \mathbf{P}\|_F^2$ simply transfers again the  user feature with a linear transformation $\mathbf{S}$. The retweeting matrix $\mathbf{X}$, similarily, also using $\mathbf{FW}^T$ as the latent vectors. 

To summarize, the model Eq.~\ref{Transform} bases the prediction of retweeting action, opinion and sentiment all on the user features. If we make $\lambda_1$ to be infinitely large, meaning that we enforce $\mathbf{FW}^T = \mathbf{O}$, then in fact, $\mathbf{X}\approx \mathbf{O}\mathbf{V}^T$ and $\mathbf{O}\mathbf{S} \approx \mathbf{P}$. Such choice is based on the assumption that opinion drives both the retweeting action and sentiment.  

Model Eq.~\ref{Transform} is an comprehensive model, in the sense that the subtask of  matrix completion, cluatering and regression are fused together, by basing all prediction on user feature transformation. What if we are not given the use feature information? Instead, we directly model the relation between retweeting action, opinion and sentiment. A straightforward model could be 

\begin{equation}
\begin{split}
    \underset{\mathbf{U,V}}{\hbox{min}} \ \mathcal{O}  & =  \|\mathbf{X}-\mathbf{U}\mathbf{V}^T\|_F^2 + \lambda_1 \|\mathbf{U} - \mathbf{O}\|_F^2 + \lambda_2 \|\mathbf{U}\mathbf{S} - \mathbf{P}\|_F^2\\& +\alpha\|\mathbf{U}\|_F^2 + \beta\|\mathbf{V}\|_F^2 +\gamma\|\mathbf{S}\|_F^2 - Tr(\Lambda_1\mathbf{U}^T) - Tr(\Lambda_2\mathbf{V}^T).
    \label{Transform_alt}
\end{split}
\end{equation}

%In this model, the relation between the three response variable is more clearly shown. Futhermore, the difference between \textit{enforcement} and \textit{transformation} is also straightforward: $\lambda_1 \|\mathbf{U} - \mathbf{O}\|_F^2$ is enforcement, while $\lambda_2 \|\mathbf{U}\mathbf{S} - \mathbf{P}\|_F^2$ is transformation.

\newpage
\bibliographystyle{plain}
\bibliography{references}
\end{document}